\documentclass[12pt]{amsart}
\usepackage{amsfonts,amssymb,amsthm,eucal,amsmath, upref, verbatim}
\usepackage[usenames, dvipsnames]{color}
\allowdisplaybreaks
\newtheorem{thm}{Theorem}
\newtheorem{cor}[thm]{Corollary}
\newtheorem{lemma}[thm]{Lemma}
\newtheorem{prop}[thm]{Proposition}

\newcommand{\R}{\mathbb{R}}
\newcommand{\C}{\mathbb{C}}

\newcommand{\E}{\mathbb{E}}

\newcommand{\inprod}[2]{\left\langle #1, #2 \right\rangle}

\newcommand{\ds}{\displaystyle}

\renewcommand{\P}{\mathbb{P}}

\renewcommand{\L}{\mathcal{L}}
\newcommand{\I}{\mathbb{I}}
\newcommand{\ee}{\mathbb{E}}

\renewcommand{\l}{\ell}
\newcommand{\ww}{W^\prime}
\newcommand{\var}{\mathrm{Var}}

\newcommand{\Hess}{\mathrm{Hess\,}}
\newcommand{\tr}{\mathrm{Tr\,}}
\newcommand{\U}{\mathcal{U}}
\renewcommand{\O}{\mathcal{O}}
\renewcommand{\Re}{\operatorname{Re}}
\renewcommand{\Im}{\operatorname{Im}}

\newcommand{\n}{\mathfrak{N}}

\begin{document}

\title{Multivariate normal approximation using exchangeable pairs}
\author{Sourav Chatterjee}\thanks{Sourav Chatterjee, 367 Evans Hall
  \#3860, Department of Statistics, University of California at
  Berkeley, Berkeley, CA 94720-3860.  Email: sourav@stat.berkeley.edu}
\author{Elizabeth Meckes}\thanks{Elizabeth Meckes, American Institute
  of Mathematics, 360 Portage Ave, Palo Alto, CA, and 220 Yost Hall,
  Department of Mathematics, Case Western Reserve University, 10900
  Euclid Ave., Cleveland, OH 44122.  Email: ese3@cwru.edu}
\begin{abstract}
Since the introduction of Stein's method in the early 1970s, much research has
been done in extending and strengthening it; however, there does not exist a 
version of Stein's original method of exchangeable pairs for multivariate normal approximation. The aim of this article is to fill this void.  We present
three abstract normal approximation theorems using exchangeable pairs
in multivariate contexts, one for situations in which the underlying
symmetries are discrete, and  real and complex
versions of a theorem for situations involving continuous 
symmetry groups.  
Our main applications are proofs of the approximate normality of rank $k$
projections of 
Haar measure on the orthogonal and unitary groups, when $k=o(n)$.
\end{abstract}
\maketitle


\section{Introduction}
Stein's method was introduced by Charles Stein \cite{stein72} as a tool for 
proving central limit theorems for sums of dependent random variables. Stein's 
version of his method, best known as the ``method of exchangeable pairs'', 
is described in detail in his later work \cite{stein86}. 
The method of exchangeable pairs is a general technique whose applicability 
is not restricted to sums of random variables; for some recent examples, one can look at the work of Jason Fulman~\cite{fulman05} on central limit 
theorems for complicated objects arising from the representation theory of 
permutation groups, and the work of the second-named author 
\cite{meckes-L} on
the distribution of eigenfunctions of the Laplacian on Riemannian manifolds.

One of the significant advantages of the method is that it automatically 
gives concrete error bounds. Although Stein's original theorem does not 
generally give Kolmogorov distance 
bounds of the correct order, there has been substantial 
research on modifications of Stein's result to obtain rate-optimal 
Berry-Ess\'een type bounds (see e.g.\ the works of Rinott \& Rotar 
\cite{rinottrotar97} and Shao \& Su \cite{shaosu04}). 
  The ``infinitesimal'' version
of the method described in \cite{meckes-thesis} and in our Theorems 
\ref{cont} and \ref{complex} below frequently does produce bounds of the
correct order, in total variation distance in the univariate case and
in Wasserstein distance in the multivariate case.  

Heuristically, the method of exchangeable pairs for univariate
normal approximation goes as follows.
Suppose that a random variable $W$ is conjectured to be approximately a
standard Gaussian. The first step in the method is to 
construct a second random variable $\ww$ on the same probability space 
such that $(W,\ww)$ is an exchangeable pair, i.e.\ $(W,\ww)$ has the 
same distribution as $(\ww,W)$. The random variable $\ww$ is generally
constructed by making a small random change in $W$, so that $W$ and $\ww$
are close. 

Let $\Delta = W-\ww$. 
The next step is to verify the existence of a small number $\lambda$ such that
\begin{align}
\ee(\Delta \mid W) &= \lambda W + r_1,\label{eq1}\\
\ee(\Delta^2\mid W) &= 2\lambda + r_2, \ \ \text{and}\label{eq2}\\
\ee|\Delta|^3 &= r_3, \label{eq3}
\end{align}
where the random quantities $r_1,r_2,$ and $r_3$ are all negligible compared 
to $\lambda$. If the above relations hold, then, depending on the sizes of 
$\lambda$ and the $r_i$'s, one can conclude that $W$ is approximately 
Gaussian. The exact statement of Stein's abstract normal approximation theorem
for piecewise differentiable test functions is the following:
\vskip.1in
\begin{thm}[Stein \cite{stein86}, page 35]\label{basic1}
Let $(W,\ww)$ be an exchangeable 
pair of real random variables such that $\ee W^2=1$ and $\ee\big[W-\ww\mid W
\big] = \lambda W$ for some $0<\lambda< 1$. Let $\Delta = W-\ww$. Let 
$h:\R\to\R$ be bounded with piecewise continuous derivative $h'$.  Then for $Z$
a standard normal random variable, 
\begin{align*}
\big|\E h(W)-\E h(Z)\big|&\le \frac{\|h-\E h(Z)\|_\infty}{\lambda}\sqrt{\var(
\ee\big[\Delta^2\mid W\big])} 
+ \frac{\|h'\|_\infty}{4\lambda} \ee|\Delta|^3.
\end{align*}
\end{thm}
\bigskip
Observe that the condition $\E\big[W-W'\big|W\big]=
\lambda W$ implies that $\E \Delta^2=2\lambda,$ and thus the bound in 
Stein's theorem above can also be stated as:
$$\big|\E h(W)-\E h(Z)\big|\le2\|h-\E h(Z)\|_\infty\sqrt{\E\left[
\frac{1}{2\lambda}\E\left[\Delta^2\big|W\right]-1\right]}+\frac{\|h'
\|_\infty}{4\lambda}\E|\Delta|^3.$$

Powerful as it is, the above theorem and all its existing
modifications cater only to {\it univariate} normal approximation.
There has been some previous work in proving multivariate central
limit theorems using Stein's method, though none of these approaches
have used exchangeable pairs.  In 1996, Rinott \& Rotar
\cite{rinottrotar96} proved multivariate central limit theorems for
sums of dependent random vectors using the dependency graph version of
Stein's method.  Around the same time, Goldstein \& Rinott
\cite{goldsteinrinott96} developed the size-bias coupling version of
Stein's method for multivariate normal approximation.  Both of these
techniques are well-known and in regular use.  More recently,
Rai{\v{c}} \cite{raic} proved a new multivariate central limit theorem
for sums of dependent random vectors with the dependency graph
approach which removed the need for finite third moments.  However, as
in the univariate case, there are many problems which are more
amenable to analysis via exchangeable pairs (particularly the
adaptation to the case of continuous symmetries) which necessitates
the creation of a multivariate version of this method. The present
authors introduced, for the first time, a multivariate version of
Theorem \ref{basic1} in an earlier draft of this manuscript that was
posted on arXiv. Subsequently, an extension of one of our main results (Theorem \ref{discrete}) to the case of
multivariate normal approximation with non-identity covariance was
formulated by Reinert and R\"ollin~\cite{reinertrollin07}. Our current
draft is mainly a reorganization of the original manuscript, with
better error bounds in several examples. Let us refer to the
Reinert-R\"ollin paper \cite{reinertrollin07} for many other
interesting applications.

The contents of this paper are as follows.  In Section \ref{results},
we prove 
three abstract normal approximation theorems which give a framework
for using the method of exchangeable pairs in a multivariate
context.  The first is
for situations in which the symmetry used in constructing the exchangeable
pair is discrete, and is a fairly direct analog of Theorem \ref{basic1} above.  
An an example, the theorem is applied in Section \ref{disc-sec} to prove a basic central limit
theorem for a sum of independent, identically distributed random vectors.

The second abstract theorem of Section \ref{results}
includes an additional modification, making
it useful in situations in which continuous symmetries are present.
The idea for the modification was introducted by Stein in \cite{steintech}
and further developed in \cite{meckes-thesis}.  Section \ref{cont-sec}
contains two applications of this theorem.  First, for $Y$ a random vector 
in $\R^n$ with 
spherically symmetric distribution, sufficient conditions are given under
which the first $k$ coordinates are approximately distributed as a standard
normal random vector in $\R^k$.  We then give a treatment of
projections of Haar measure on the orthogonal group.  
Specifically, for $M$ a random $n\times n$
orthogonal matrix and $A_1,\ldots,A_k$ fixed matrices over $\R$, 
we give an explicit bound on the Wasserstein distance between 
$(\tr(A_1M),\ldots,\tr(A_kM))$ and a Gaussian random vector.  

As a corollary to the theorem discussed above, we state a theorem for
bounding the distance between a complex random vector and a complex Gaussian
random vector, in the context of continuous groups of symmetries.  The
main application of this version of the theorem in given in Section 
\ref{cont-sec}, where for $M$ a random $n\times n$ unitary 
matrix and $A_1,\ldots,A_n$ fixed matrices over $\C$, we derive an explicit
bound on the Wasserstein distance between $(\tr(A_1M),\ldots,\tr(A_kM))$ 
and a complex Gaussian random vector.  

\bigskip

Before moving into Section \ref{results}, we give the following very brief 
outline of the literature around the various other versions of Stein's method.
\vskip.1in
\noindent {\bf Other versions of Stein's method.}
The three most notable variants of Stein's method are (i) the dependency 
graph approach introduced by Baldi and Rinott \cite{baldirinott89} and 
 further developed by Arratia, Goldstein and 
Gordon \cite{agg90} and
Barbour, Karo{\'n}ski, and Ruci{\'n}ski \cite{bkr89},
(ii) the size-biased coupling method of Goldstein and Rinott \cite{goldsteinrinott97} (see also Barbour, Holst and Janson \cite{bhj92}), and (iii) the zero-biased coupling technique due to 
Goldstein and Reinert \cite{goldsteinreinert97}.  In addition to these three
basic approaches, an important contribution was made by 
Andrew Barbour \cite{barbour90}, who noticed the connection between 
Stein's method and diffusion approximation. This connection has 
subsequently been widely exploited by practitioners of Stein's method, and is 
a mainstay of some of our proofs.

Besides normal approximation, Stein's method has been successfully used for 
proving convergence to several other distributions as well. 
Shortly after
the method was introduced for normal approximation by Stein, Poisson 
approximation by Stein's method was introduced by Chen \cite{chen75} and 
became popular after the publication of \cite{agg89, agg90}.  
The method has also been developed for gamma 
approximation by Luk \cite{luk94}; for chi-square approximation by 
Pickett \cite{pick}; for the uniform distribution on the discrete circle 
by Diaconis \cite{diaconis04}; for the semi-circle law by G\"otze and Tikhomirov
\cite{gottik05}; for the binomial and multinomial distributions by Holmes
\cite{holmes04} and Loh \cite{loh92}; and the hypergeometric distribution, also
by Holmes \cite{holmes04}.

The method of exchangeable pairs was extended to Poisson approximation by Chatterjee, Diaconis and Meckes in the survey paper \cite{cdm05}, and to a general method of normal approximation for arbitrary functions of independent random variables in \cite{chatterjee06}.

For further references and exposition (particularly to the method of exchangeable pairs), we refer to the recent monograph \cite{diaconisholmes04}. 

\subsection{Notation and conventions}

The total variation distance $d_{TV}(\mu,\nu)$ between the measures $\mu$ and
$\nu$ on $\R$ is defined by
$$d_{TV}(\mu,\nu)=\sup_A\big|\mu(A)-\nu(A)|,$$
where the supremum is over measurable sets $A$.  This is equivalent to
$$d_{TV}(\mu,\nu)=\frac{1}{2}\sup_{f}\left|\int f(t)d\mu(t)-
\int f(t)d\nu(t)\right|,$$
where the supremum is taken over continuous functions which are bounded
by 1 and vanish at infinity; this is the definition most commonly used 
in what follows. 
The total variation distance between two random variables $X$ and $Y$ is
defined to be the total variation distance between their distributions:
$$d_{TV}(X,Y)=\sup_A\big|\P(X\in A)-\P(Y\in A)\big|=\frac{1}{2}\sup_{f}
\big|\E f(X)-\E f(Y)\big|.$$
If the Banach space of signed measures on $\R$ is viewed as dual to the space
of continuous functions on $\R$ vanishing at infinity, then the total variation
 distance is (up to the factor of $\frac{1}{2}$) the norm distance on that 
Banach space.

The Wasserstein distance 
$d_{W}(X,Y)$ between the random variables $X$ and $Y$ is defined by
$$d_{W}(X,Y)=\sup_{M_1(g)\le1}\big|\E g(X)-\E g(Y)\big|,$$
where $M_1(g)=\sup_{x\neq y}\frac{|g(x)-g(y)|}{|x-y|}$ is the Lipschitz 
constant of $g$.  Note that Wasserstein distance is not directly comparable
to total variation distance, since the class of functions considered is 
required to be Lipschitz but not required to be bounded.   In particular, total
variation distance is always bounded by 1, whereas the statement that the 
Wasserstein distance between two distributions is bounded by 1 has content.
  On the space of 
probability distributions
with finite absolute first moment, Wasserstein
distance induces a stronger topology than the usual 
one described by weak convergence, but not as strong as 
the topology induced by the total variation distance.  See \cite{dud} for 
detailed discussion of the various notions of distance between probability
distributions.

We will use $\n(\mu,\sigma^2)$ to denote the normal distribution on $\R$ with 
mean $\mu$ and variance $\sigma^2$; unless otherwise stated, the 
random variable $Z=(Z_1,\ldots,Z_k)$ is  
understood to be a standard Gaussian random vector on $\R^k$.  

In $\R^n$, the Euclidean inner product is denoted $\inprod{\cdot}{\cdot}$ and
the Euclidean norm is denoted $|\cdot|$.  
On the space of real (resp. complex) $n\times n$ matrices, 
the Hilbert-Schmidt inner
product is defined by
$$\inprod{A}{B}_{H.S.}=\tr(AB^T),\qquad\qquad\big({\rm resp.}\,\inprod{A}{
B}_{H.S.}=\tr(AB^*)\big)$$
with corresponding norms
$$\|A\|_{H.S.}=\sqrt{\tr(AA^T)},\qquad\qquad\left({\rm resp.}\,\|A\|_{H.S.}=\sqrt{
\tr(AA^*)}\right).$$
The operator norm of a matrix $A$ over $\R$ is defined by
$$\|A\|_{op}=\sup_{|v|=1,|w|=1}|\inprod{Av}{w}|.$$
The $n\times n$ identity matrix is denoted $I_n$, the 
$n\times n$ matrix of all zeros is denoted $0_n$, and 
$A\oplus B$ is the block direct
sum of $A$ and $B$.

For $\Omega$ a domain in $\R^k$, the notation 
$C^k(\Omega)$ will be used for the space of $k$-times 
continuously differentiable real-valued functions on $\Omega$, and 
$C^k_o(\Omega)\subseteq C^k(\Omega)$ are those $C^k$ functions on 
$\Omega$ with compact support.
For $g:\R^k\to\R$, let
$$M_1(g):=\sup_{x\neq y}\frac{|g(x)-g(y)|}{|x-y|};$$
if $g\in\C^1(\R^k)$ also, then let 
$$M_2(g):=\sup_{x\neq y}\frac{|\nabla g(x)-\nabla g(y)|}{|x-y|};$$
if $g\in\C^2(\R^k)$ as well, then
$$M_3(g):=\sup_{x\neq y}\frac{\|\Hess g(x)-\Hess g(y)\|_{op}}{|x-y|}.$$

The last definition differs from the one in \cite{raic}, where $M_3$
is defined in terms of the Hilbert-Schmidt norm as opposed to the 
operator norm.  
Note that if $g\in C^1(\R^k)$, then $M_1(g)=\sup_x|\nabla g(x)|$, and if 
$g\in C^2(\R^k)$, then $M_2(g)=\sup_x\|\Hess g(x)\|_{op}.$

\medskip


\section{Two abstract normal approximation theorems}\label{results}
In this section we develop the general machine that  will be applied in
the examples in Sections \ref{disc-sec} and \ref{cont-sec}. In the following, we use the notation $\L (X)$ to denote the law of a random vector or variable $X$.
 The following lemma gives a second-order
characterizing operator for the Gaussian distribution on $\R^k$. 
\begin{lemma}\label{char}
Let $Z\in\R^k$ be a random vector with $\{Z_i\}_{i=1}^k$ independent, identically
distributed standard Gaussians.  
\begin{enumerate}
\item \label{char1}
If $f:\R^k\to\R$ is two times continuously differentiable and 
compactly supported, then 
$$\E\big[\Delta f(Z)-\inprod{Z}{\nabla f(Z)}\big]=0.$$

\item \label{char2}If $Y\in\R^k$ is a random vector such that 
$$\E\big[\Delta f(Y)-\inprod{Y}{\nabla f(Y)}\big]=0$$
for every $f\in C^2(\R^k)$ with $\E\big|\Delta f(Y)-\inprod{Y}{\nabla f(Y)}
\big|<\infty$, then $\L(Y)=\L(Z)$.  

\item \label{sol}If $g\in C_o^\infty(\R^k)$, then the function
$$U_og(x):=\int_0^1\frac{1}{2t}\big[\E g(\sqrt{t}x+\sqrt{1-t}Z)-\E g(Z)\big]
dt$$
is a solution to the differential equation
\begin{equation}\label{diffeq}
\Delta h(x)-\inprod{x}{\nabla h(x)}=g(x)-\E g(Z).\end{equation}
\end{enumerate}
\end{lemma}
\noindent {\it Remark.} The form of $U_o g$ is a direct rewriting of the inverse of the Ornstein-Uhlenbeck generator (see Barbour \cite{barbour90}). 
\begin{proof}
Part \ref{char1} is just integration by parts.  

Part (ii) follows easily from part (iii): note that if 
$$\E\big[\Delta f(Y)-\inprod{Y}{\nabla 
f(Y)}\big]=0$$
for every $f\in C^2(\R^k)$ with $\E\big|\Delta f(Y)-\inprod{Y}{\nabla 
f(Y)}\big|<\infty$, then for $g\in C_o^\infty$ given,
$$\E g(Y)-\E g(Z)=\E\big[\Delta (U_og)(Y)-\inprod{Y}{\nabla (U_og)(Y)}\big]=0,$$
and so $\L(Y)=\L(Z)$ since $C_o^\infty$ is dense in the class of bounded
continuous functions vanishing at infinity, with respect to the supremum norm.

A proof of part (iii) is given in \cite{barbour90}, \cite{got91} and 
\cite{raic}, all using results about Markov semi-groups.  
For a direct proof, see 
\cite{meckes-thesis}.  

\end{proof}

The next lemma gives useful bounds on $U_og$ and its derivatives
in terms of $g$ and its derivatives.  
As in \cite{raic},  bounds are most 
naturally given
in terms of the 
quantities $M_i(g)$ defined in the introduction.
\begin{lemma}\label{bounds}
For $g:\R^k\to\R$ given, $U_og$ satisfies the following bounds:
\begin{enumerate}
\item \label{Hessbd}$$\sup_{x\in R^k}\|\Hess U_og(x)\|_{H.S.}\le M_1(g).$$
\item \label{M3bd}$$M_3(U_og)\le\frac{\sqrt{2\pi}}{4}M_2(g).$$

\end{enumerate}
\end{lemma}

\begin{proof}
Write $h(x)=U_og(x)$ and $Z_{x,t}=\sqrt{t}x+\sqrt{1-t}Z$.  
Note that by the formula for $U_og$,
\begin{equation}\label{derivs}
\frac{\partial^rh}{\partial x_{i_1}\cdots\partial x_{i_r}}(x)=
\int_0^1(2t)^{-1}t^{r/2}\E\left[\frac{\partial^rg}{\partial x_{i_1}
\cdots\partial x_{i_r}}(Z_{x,t})\right]dt.
\end{equation} 
It follows by integration by parts on the Gaussian expectation that
\begin{equation}\begin{split}\label{seconds}
\frac{\partial^2h}{\partial x_i\partial x_j}(x)&=\int_0^1\frac{1}{2}\E\left[
\frac{\partial^2g}{\partial x_i\partial x_j}(\sqrt{t}x+\sqrt{1-t}Z)\right]dt\\
&=\int_0^1\frac{1}{2\sqrt{1-t}}\E\left[Z_i\frac{\partial g}{\partial x_j}(Z_{x,t})
\right]dt,
\end{split}\end{equation}
and so
\begin{equation}\label{Hessian}
\Hess h(x)=\int_0^1\frac{1}{2\sqrt{1-t}}\E\left[Z\left(\nabla g(
Z_{x,t})\right)^T\right]dt.\end{equation}

Fix a $k\times k$ matrix $A$.  Then
\begin{equation*}\begin{split}
\inprod{\Hess h(x)}{A}_{H.S.}&=\int_0^1\frac{1}{2\sqrt{1-t}}\E\left[\inprod{
A^TZ}{\nabla g(Z_{x,t})}\right]dt,
\end{split}\end{equation*}
thus
$$\left|\inprod{\Hess h(x)}{A}_{H.S.}\right|\le
M_1(g)\E|A^TZ|\int_0^1\frac{1}{2\sqrt{1-t}}dt=M_1(g)\E|A^TZ|.$$
If $A=[a_{ij}]_{i,j=1}^k,$ then
$$\E|AZ|\le\sqrt{\E|AZ|^2}=\sqrt{\E\sum_{i=1}^k\left(\sum_{j=1}^ka_{ji}Z_j
\right)^2}=\sqrt{\sum_{i,j=1}^ka_{ij}^2}=\|A\|_{H.S.},$$
and thus
$$\|\Hess h(x)\|_{H.S.}\le M_1(g)$$
for all $x\in\R^k$, hence part \ref{Hessbd}.

For part \ref{M3bd}, let $u$ and $v$ be fixed vectors in $\R^k$ 
with $|u|=|v|=1.$  Then it follows from \eqref{Hessian} that 
$$\inprod{\left(\Hess h(x)-\Hess h(y)\right)u}{v}=\int_0^1\frac{1}{2\sqrt{1-t}
}\E\left[\inprod{Z}{v}\inprod{\nabla g(Z_{x,t})-\nabla g(Z_{y,t})}{u}\right]
dt,$$
and so
\begin{equation*}\begin{split}
\left|\inprod{(\Hess h(x)-\Hess h(y))u}{v}\right|&\le |x-y|\,M_2(g)\,\E|\inprod{Z}{v}|
\int_0^1\frac{\sqrt{t}}{2\sqrt{1-t}}dt\\&=|x-y|\,M_2(g)\frac{\sqrt{2\pi}}{4},
\end{split}\end{equation*}
since $\inprod{Z}{v}$ is just a standard Gaussian random variable.  This 
completes the proof of part \ref{M3bd}.
\end{proof}

There is an important difference in the behavior of solutions to the Stein
equation \ref{sol} in the context of multivariate approximation versus 
univariate approximation.  In the univariate case, one can replace the 
expression on the left-hand side of \ref{sol} with the first-order
expression
$h'(x)-xh(x);$
the function $g(x)=U_oh(x)$ which solves the differential equation
$$h'(x)-xh(x)=g(x)-\E g(Z)$$
satisfies the bounds (see \cite{stein86})
$$\|g\|_\infty\le\sqrt{\frac{\pi}{2}}\|h-\E h(Z)\|_\infty\qquad
M_1(g)\le2\|h-\E h(Z)\|_\infty\qquad M_2(g)\le2M_1(h),$$
and the fact that the differential equation is first order rather than second
then allows for reducing the degree of smoothness needed by one, over what is
required in the multivariate case.  Alternatively, one can use the same 
expression as in \ref{sol} above; in this case, $M_3(g)\le 2M_1(g)$ (see 
\cite{raic}), also decreasing by one the degree of smoothenss needed.  This
improvement allowed the univariate version \cite{M}
of Theorem \ref{ind} below, on the approximation of projections of Haar measure
on the orthogonal group by Gaussian measure, to 
be proved in total variation distance as opposed
to Wasserstein distance.

This improvement is not possible in the multivariate case; it can be shown,
for example (see \cite{raic}), that if $$f(x,y)=\max\{\min\{x,y\},0\},$$
then $U_of$ defined as in Lemma \ref{char} is twice differentiable but
$\frac{\partial^2(U_of)}{\partial x^2}$ is not Lipschitz.

\bigskip

\begin{thm}\label{discrete}
Let $X$ and $X'$ be two random vectors in $\R^k$ such that $\L(X)=\L(X')$, 
and let $Z=(Z_1,\ldots,Z_k)\in\R^k$ be a standard Gaussian random vector.  
Suppose there is a constant $\lambda$ such that 
\begin{equation}\frac{1}{\lambda}\E\left[X'-X\big|X\right]=-X\label{lin-cond}.
\end{equation}
Define the random matrix $E$ by 
\begin{equation}\frac{1}{2\lambda}\E\left[(X'-X)(X'-X)^T\big|X\right]=\sigma^2
I_k+\E\left[E\big|X\right].\label{quad-cond}\end{equation}

Then if $g\in C^2(\R^k)$ with $M_1(g)<\infty$ and $M_2(g)<\infty$,  
\begin{equation}\begin{split}\label{discbound}
\big|\E g(X) -\E g(\sigma Z)\big| &\le\frac{1}{\sigma}M_1(g)\E\|E\|_{H.S.}+
\left(\frac{\sqrt{2\pi}}{24\sigma}\right)\frac{M_2(g)}{\lambda}\E|X'-X|^3.
\end{split}\end{equation}

\end{thm}
\begin{proof}Fix $g$, and let
$U_og$ be as in Lemma \ref{char}.  Note that it suffices to assume that
$g\in C^\infty(\R^k)$:  let $h:\R^k\to\R$ be a centered Gaussian density with
covariance matrix $\epsilon^2I_k$.  Approximate $g$ by 
$g*h$; clearly  $\|g*h-g\|_\infty\to0$ as $\epsilon
\to0$, and by Young's inequality, $M_1(g*h)\le M_1(g)$ and 
$M_2(g*h)\le M_2(g)$. 

Note also that if $f(x)=g(\sigma x)$, then
$\big|\E g(X)-\E g(\sigma Z)\big|=\big|\E f(\sigma^{-1}X)-\E f(Z)\big|.$
It is easy to see that $M_1(f)=\sigma M_1(g)$ and $M_2(f)=\sigma^2 M_2(g)$.
It thus follows from the theorem for $\sigma=1$
that 
\begin{equation*}\begin{split}
\big|\E g(X)-\E g(\sigma Z)\big|&\le\sigma M_1(g)\E\left\|\sigma^{-2}E
\right\|_{H.S.}+\left(\frac{\sqrt{2\pi}}{24}\right)\frac{\sigma^2M_2(g)}
{\lambda}\E\left|\sigma^{-3}(X'-
X)\right|^3\\&=\frac{M_1(g)}{\sigma}\E\|E\|_{H.S.}+
\frac{\sqrt{2\pi}\,M_2(g)}{24\sigma\lambda}\E|X'-X|^3;
\end{split}\end{equation*}
we therefore restrict our attention to the case $\sigma=1$.

For notational convenience, write 
$h(x)=U_og(x)$.  Then
\begin{eqnarray}
0&=&\frac{1}{\lambda}\E\left[h(X')-h(X)\right]\nonumber\\
&=&\frac{1}{\lambda}\E\left[\inprod{X'-X}{\nabla h(X)}+\frac{1}{2}(X'-X)^T(
\Hess h(X))(X'-X)+R\right]\nonumber\\
&=&\E\left[-\inprod{X}{\nabla h(X)}+\Delta h(X) +\inprod{E}{\Hess h(X)}_{H.S.} 
+\frac{R}{\lambda}\right]\nonumber\\
&=&\E g(X)-\E g(Z)+\E\left[\inprod{E}{\Hess h(X)}_{H.S.} 
+\frac{R}{\lambda}\right], \label{diff}
\end{eqnarray}
where $R$ is the error in the second-order expansion.  
By an alternate form of Taylor's theorem (see \cite{yom83}),
\begin{equation*}\begin{split}
\E|R|\le\frac{M_3(h)}{6}\E|X'-X|^3\le\frac{\sqrt{2\pi}M_2(g)}{24}\E|X'-X|^3.
\end{split}\end{equation*}
Furthermore, 
$$\E\left|\inprod{E}{\Hess h(X)}\right|\le \left(\sup_{y\in \R^k}\|\Hess h(y)
\|_{H.S.}\right)
\E\|E\|_{H.S.}\le M_1(g)\E\|E\|_{H.S.}.$$
This completes the proof.

\end{proof}

\medskip

\noindent{\em Remarks.}
\begin{enumerate}
\item Usually the $X$ and $X'$ of the theorem will make an exchangeable
pair, but this is not required for the proof.  
\item The coupling assumed in \eqref{lin-cond} implies that $\E X=0.$  It is
not required that $X$ have a scalar covariance matrix, however, it follows
from \eqref{lin-cond} and \eqref{quad-cond} that 
$$\E \big[E\big]=\E\big[XX^T\big]-\sigma^2I_k.$$
It should therefore be the case that the covariance matrix of $X$ is not 
too far from $\sigma^2I_k$.  
\end{enumerate}

\medskip

The following is a continuous analog of Theorem \ref{discrete}.  A univariate
version which gives approximation in total variation distance was
proved in \cite{M}.  
As was noted following the proof of Lemma \ref{bounds},
a bound on total variation distance in the multivariate context is not 
possible with the method used here because of the difference in the behavior
of solutions to the Stein equation in the multivariate context.

\begin{thm}\label{cont}
Let $X$ be a random vector in $\R^k$ and for each $\epsilon >0$ let $X_\epsilon$ be a random vector such that $\L(X)=
\L(X_\epsilon)$, with the property that $\lim_{\epsilon\to0}X_\epsilon=X$ 
almost surely.  Let $Z$ be a normal random vector in $\R^k$ with mean zero and
covariance matrix $\sigma^2I_k$.  
Suppose there is a function $\lambda(\epsilon)$ and a random matrix
$F$ such that the following conditions hold.
\begin{enumerate}
\item $$\frac{1}{\lambda(\epsilon)}\E\left[(X_\epsilon-X)_i\big|X\right]
\xrightarrow[\epsilon\to0]{L_1}- X. $$\label{first-diff}
\item $$\frac{1}{2\lambda(\epsilon)}\E\left[(X_\epsilon-X)(X_\epsilon-X)^T|
X\right]\xrightarrow[\epsilon\to0]{L_1}\sigma^2I_k+\E\left[F\big|X\right].$$
\label{second-diff}
\item For each $\rho>0$, $$\lim_{\epsilon\to0}\frac{1}{\lambda(\epsilon)}
\E\left[\big|X_\epsilon-X\big|^2
\I(|X_\epsilon-X|^2>\rho)\right]=0.$$\label{lindeberg}
\end{enumerate}
Then
\begin{equation}\label{contbd}
d_W(X,Z)\le\frac{1}{\sigma}\E\|F\|_{H.S.}
\end{equation}
\end{thm}

\begin{proof}Fix a test function $g$; as in the proof of Theorem 
\ref{discrete}, it suffices to assume that $g\in C^\infty(\R^k)$ and to
consider only the case $\sigma=1$; the general result follows exactly as
before.

Let $U_og$ be as in Lemma \ref{char}, and as before, write
$h(x)=U_og(x)$.  Observe
\begin{equation}\label{earlyexp}
\begin{split}
0&=\frac{1}{\lambda(\epsilon)}\E\left[h(X_\epsilon)-h(X)\right]\\
&=\frac{1}{\lambda(\epsilon)}
\E\left[\inprod{X_\epsilon-X}{\nabla h(X)}+\frac{1}{2}(X_
\epsilon-X)^T(\Hess h(X))(X_\epsilon-X)+R\right],
\end{split}
\end{equation}
where $R$ is the error in the second-order approximation of 
$h(X_\epsilon)-h(X)$.  By Taylor's theorem, there is a constant $K$ 
(depending on $h$) and
a function $\delta$ with $\delta(x)\le K\min\{x^2,x^3\},$ such that  
$\big|R\big|\le \delta(|X'-X|).$
Fix $\rho>0$.  Then by breaking up the integrand over the sets
$\big\{|X_\epsilon-X|
\le\rho\big\}$ and $\big\{|X_\epsilon-X|>\rho\big\}$, 
\begin{equation*}\begin{split}
\frac{1}{\lambda(\epsilon)}\E\big|R\big|&\le\frac{K}{\lambda(\epsilon)}
\E\Big[|X_\epsilon-X|^3\I(|X_\epsilon-X|\le\rho)+|X_\epsilon-X|^2\I(|X_\epsilon
-X|>\rho)\Big]\\&\le \frac{K\rho\E\big|X_\epsilon-X\big|^2}{\lambda(\epsilon)}
+\frac{K}{\lambda(\epsilon)}\E\Big[|X_\epsilon-X|^2\I(|X'-X|>\rho)\Big].
\end{split}\end{equation*}
The second term tends to zero as $\epsilon\to0$ by condition \ref{lindeberg};
condition \ref{second-diff} implies that the first is bounded by $CK\rho$ for 
a constant $C$ depending on $k$ and on the 
distribution of $X$.  It follows that 
$$\lim_{\epsilon\to0}\frac{1}{\lambda(\epsilon)}\E\big|R\big|=0.$$
For the first two terms of \eqref{earlyexp},
\begin{equation}\begin{split}
\lim_{\epsilon\to0}&\frac{1}{\lambda(\epsilon)}\E\left[\inprod{X_\epsilon-X}{
\nabla h(X)}+\frac{1}{2}(X_\epsilon-X)^T(\Hess h(X))(X_\epsilon-X)\right]\\&=
\lim_{\epsilon\to0}\frac{1}{
\lambda(\epsilon)}\E\left[\inprod{\E\left[(X_\epsilon-X)\big|X\right]}{\nabla
h(X)}+\frac{1}{2}\inprod{\E\left[\big(X_\epsilon-X\big)\big(X_\epsilon-X\big)^T
\big|X\right]}{\Hess h(X)}_{H.S.}\right]\\
&=\E\left[-\inprod{X}{\nabla h(X)}+\Delta h(X) +
\inprod{\E\left[F\big|X\right]}{\Hess h(X)}_{H.S.}\right]\\
&=\E g(X)-\E g(Z)+\E\left[\inprod{\E\left[F\big|X\right]}
{\Hess h(X)}_{H.S.}\right], 
\end{split}\end{equation}
where
conditions \ref{first-diff} and \ref{second-diff} together with the boundedness
of $\nabla h$ and $\Hess h$ are used to get the third line and the definition
of $h=U_og$ is used to get the fourth line.
We have thus shown that 
\begin{equation}
\E\big[g(X)-g(Z)\big]=-\E\inprod{F}{\Hess h(X)}_{H.S.}.
\label{equality}
\end{equation}
The result now follows immediately by applying the Cauchy-Schwarz inequality
to \eqref{equality} and then the bound $\|\Hess h(x)\|_{H.S.}\le M_1(g)$ from
Lemma \ref{bounds} \ref{Hessbd}.

\end{proof}

\noindent{\em Remarks.}
\begin{enumerate}
\item It is easy to see that if
$$\begin{array}{ll}{\rm (iii')}\qquad
 \lim_{\epsilon\to0}\frac{1}{\lambda(\epsilon)}\E\big|X_\epsilon-X
\big|^3=0,&\phantom{{\rm (iii')}\qquad
 \lim_{\epsilon\to0}\frac{1}{\lambda(\epsilon)}\E\big|X_\epsilon-X
\big|^3=0}\end{array}$$
then condition (iii) of the theorem holds.  This is what is done in the 
applications below.

\item As in Theorem \ref{discrete}, the condition (i) implies that $\E X=0$ and
it follows from (i) and (ii) that 
$$\E F=\E XX^T-\sigma^2I;$$
the covariance matrix of $X$ should thus not be far from $\sigma^2I$.

\end{enumerate}

\medskip

Theorem \ref{cont} has the following corollary for complex random vectors.

\begin{cor}\label{complex}

Let $W$ be a random vector in $\C^k$ and for each $\epsilon >0$ let 
$W_\epsilon$ be a random vector such that $\L(W)=
\L(W_\epsilon)$, with the property that $\lim_{\epsilon\to0}W_\epsilon=W$ 
almost surely.  Let $Z=(Z_1,\ldots,Z_k)$ be a 
standard complex Gaussian random vector; i.e., 
with covariance matrix of
the corresponding random vector in $\R^{2k}$ given by
$\frac{1}{2}I_{2k}$.  Suppose there is a function $\lambda(\epsilon)$ and 
complex $k\times k$ random matrices $\Gamma=\big[\gamma_{ij}\big]$ and 
$\Lambda=\big[\lambda_{ij}\big]$ such that
\begin{enumerate}
\item $$\frac{1}{\lambda(\epsilon)}\E\left[(W_\epsilon-W)\big|W\right]
\xrightarrow[\epsilon\to0]{L_1}- W. $$\label{first-diff-c}
\item $$\frac{1}{2\lambda(\epsilon)}\E\left[(W_\epsilon-W)(W_\epsilon-W)^*|
W\right]\xrightarrow[\epsilon\to0]{L_1}I_k+\E\left[\Gamma\big|W\right].$$
\label{second-diff-c-1}

\item $$\frac{1}{2\lambda(\epsilon)}\E\left[(W_\epsilon-W)(W_\epsilon-W)^T|
W\right]\xrightarrow[\epsilon\to0]{L_1}\E\left[\Lambda\big|W\right].$$
\label{second-diff-c-2}

\item For each $\rho>0$, $$\lim_{\epsilon\to0}\frac{1}{\lambda(\epsilon)}
\E\left[\big|W_\epsilon-W\big|^2
\I(|W_\epsilon-W|^2>\rho)\right]=0.$$
\end{enumerate}
Then 
$$d_W(W,Z)\le\E\|\Gamma\|_{H.S.}+\E\|\Lambda\|_{H.S.}.$$
\end{cor}

\begin{proof}
Identifying $\C^k$ with $\R^{2k}$, $W$ satisfies the conditions of Theorem 
\ref{cont} with $\sigma^2=\frac{1}{2}$ and $F$ given as a $k\times k$ matrix
of $2\times 2$ blocks, with the $i$-$jth$ block equal to
$$\frac{1}{2}\begin{bmatrix}\Re(\gamma_{ij}+\lambda_{ij}) & \Im(\lambda_{ij}-
\gamma_{ij})\\\Im(\lambda_{ij}+\gamma_{ij}) 
& \Re(\gamma_{ij}-\lambda_{ij})\end{bmatrix}.$$
Thus $\|F\|_{H.S.}^2=\frac{1}{2}(\|\Gamma\|_{H.S.}^2+\|\Lambda\|_{H.S.}^2)$ and 
$$\E\|F\|_{H.S.}\le\frac{1}{\sqrt{2}}\Big[\E\|\Gamma\|_{H.S.}+\E\|\Lambda\|_{H.S.}\Big].$$
\end{proof}

\section{Examples using Theorem \ref{discrete}}\label{disc-sec}
\subsection{A basic central limit theorem}
As a simple illustration of the use of Theorem \ref{discrete}, we
derive error bounds in the classical multivariate CLT for sums of
independent random vectors. While the question of error bounds in the
univariate CLT was settled long ago, the optimal bounds in the
multivariate case are still unknown 
and much work has been done in this
direction.  One important contribution was made by G\"otze \cite{got91},
who used Stein's method in conjunction with induction.  
To the best of our knowledge, the
most recent results are due to V.~Bentkus \cite{bentkus03}, where one
can also find extensive pointers to the literature.

Suppose $Y$ is a random vector in $\R^k$ with mean zero and identity
covariance. Let $W$ be the normalized sum of $n$ i.i.d.\ copies of
$Y$. G\"otze \cite{got91} and Bentkus \cite{bentkus03} both give bounds on 
quantities like 
$\Delta_n = \sup_{f\in \mathcal{A}} |\E f(W)-\E f(Z)|$,
where $Z=(Z_1,\ldots,Z_k)$ is a standard $k$-dimensional normal random 
vector and
$\mathcal{A}$ is any collection of functions satisfying certain
properties. For example, when $\mathcal{A}$ is the class of indicator 
functions of convex
sets, Bentkus gets $\Delta_n \le 400k^{1/4}n^{-1/2}\ee|Y|^3$, improving
on G\"otze's earlier bound which has a coefficient of $k^{1/2}$ rather
than $k^{1/4}$. Note
that $\ee|Y|^3 = O(k^{3/2})$.

Theorem \ref{discrete} allows us to easily obtain uniform bounds on $|\E g(S_n) - \ee g(Z)|$ for large classes of smooth functions.
\begin{thm}\label{basic}
Let $\{Y_i\}_{i=1}^n$ be a set of independent, identically distributed
 random vectors in $\R^k$.   Assume
that the $Y_i$ are such that 
$$\E (Y_i)=0,\qquad\E (Y_iY_i^T)=I_k.$$
  Let $W=\frac{1}{\sqrt{n}}\sum_{i=1}^nY_i.$  Then for any $g\in C_o^2$, 
$$\big|\E g(W)-\E g(Z)\big| \le \frac{M_1(g)}{2\sqrt{n}}\sqrt{\E|Y_1|^4-k}+
\frac{\sqrt{2\pi}}{3\sqrt{n}}M_2(g)\E|Y_1|^3.$$

\end{thm}
\begin{proof}
To apply Theorem \ref{discrete}, make an exchangeable pair $(W,W')$ as 
follows.  For each $i$, let $X_i$ be an independent copy of $Y_i$, and 
let $I$ be a uniform random variable in $\{1,\ldots,n\}$, 
independent of everything. Define $W'$ by
$$W'=W-\frac{Y_I}{\sqrt{n}}+\frac{X_I}{\sqrt{n}}.$$
Then
\begin{equation*}\begin{split}
&\E\left[W'-W\big|W\right]=\frac{1}{\sqrt{n}}\E\left[X_I-Y_I\big|W\right]\\
&=\frac{1}{n^{3/2}}\sum_{i=1}^n\E\left[X_i-Y_i\big|W\right] =-\frac{1}{n}W,
\end{split}\end{equation*}
where the independence of $X_i$ and $W$ has been used in the last line.
Thus condition \ref{lin-cond} of Theorem \ref{discrete} holds with $\lambda=\frac{1}{n}$.

It remains to check condition 2 and bound the $E_{ij}$.  
Write $Y_i=(Y_i^1,\ldots,Y_i^k)$.  For $1\le j,\ell\le k $,
\begin{equation*}\begin{split}
E_{j\ell}&=\frac{n}{2}\E\left[(W_j'-W_j)(W_\ell'-W_\ell)\big|W\right]
-\delta_{j\ell}\\
&=\frac{1}{2}\E\left[(X^j_I-Y^j_I)(X^\ell_I-Y^\ell_I)\big|W\right]-\delta_{j\ell}
\\&=\frac{1}{2n}\sum_{i=1}^n\E\left[X_i^jX_i^\ell-X_i^jY_i^\ell-X_i^\ell 
Y_i^j+Y_i^jY_i^\ell\big|W\right]-\delta_{j\ell}\\
&=\frac{1}{2n}\sum_{i=1}^n\E\left[Y_i^jY_i^\ell-\delta_{j\ell}\big|W\right],
\end{split}\end{equation*}
by the independence of the $X_i$ and the $Y_i$.  
Thus
\begin{equation*}\begin{split}
\E E_{j\ell}^2&=\frac{1}{4n^2}\E\left(\E\left[\left.\sum_{i=1}^n(Y_i^jY_i^\ell
-\delta_{j\ell})\right|W\right]\right)^2\\&\le\frac{1}{4n^2}\E\left[\left(
\sum_{i=1}^n(Y_i^jY_i^\ell-\delta_{j\ell})\right)^2\right]\\&=\frac{1}{4n^2}
\E\sum_{i=1}^n(Y_i^jY_i^\ell-\delta_{j\ell})^2\\&=\frac{1}{4n}\E\left[Y_1^jY_1^\ell-\delta_{
j\ell}\right]^2\\&=\frac{1}{4n}\left[\E\left(Y_1^jY_1^\ell\right)^2-\delta_{j
\ell}\right],
\end{split}\end{equation*}
where the independence of the $Y_i$ has been used to get the third line.
It follows that 
\begin{equation*}\begin{split}
\E\|E\|_{H.S.}\le\sqrt{\E\|E\|_{H.S.}^2}&\le\frac{1}{2\sqrt{n}}\sqrt{\sum_{j,
\ell}\left(\E(Y_1^jY_1^\ell)^2-\delta_{j\ell}\right)}=\frac{1}{2\sqrt{n}}
\sqrt{\E|Y_1|^4-k}.
\end{split}\end{equation*}

It remains to bound the second term of Theorem \ref{discrete}.
\begin{equation*}\begin{split}
\frac{1}{\lambda}\E|W'-W|^3&=\frac{1}{\sqrt{n}}\E\big|X_I-Y_I\big|^3\\&=
\frac{1}{\sqrt{n}}\E\big|X_1-Y_1\big|^3\\&\le\frac{1}{\sqrt{n}}\E\left(
|X_1|^3+3|X_1|^2|Y_1|+3|Y_1|^2|X_1|+|Y_1|^3\right).
\end{split}\end{equation*}
Applying H\"older's inequality with $p=\frac{3}{2}$ and $q=3$, 
$$\E|X_1|^2|Y_1|\le\left(\E|X_1|^3\right)^{2/3}\left(\E|Y_1|^3\right)^{1/3}
=\E|Y_1|^3.$$
It follows that 
$$\frac{1}{\lambda}\E|W'-W|^3\le\frac{8\E|Y_1|^3}{\sqrt{n}}.$$
Together with Theorem \ref{discrete}, this finishes the proof.
\end{proof}

\section{Examples using Theorem \ref{cont}}\label{cont-sec}

\subsection{Rank $k$ projections of spherically symmetric measures on $\R^n$}
\label{sphere}
Consider a random vector $Y\in\R^{n}$ whose distribution is spherically
symmetric; i.e., if $U$ is a fixed orthogonal matrix, then the distribution
of $Y$ is the same as the distribution of $UY$.  
Assume that $Y$ is 
normalized such that $\E Y_1^2=1$.  Note 
that the spherical symmetry then implies that $\E YY^T=I_n$.
Assume further that there is a constant $a$ (independent of $n$) so that
\begin{equation}\label{varbda}\var(|Y|^2)\le a.\end{equation}
For $k$ fixed, let $P_k$ denote the orthogonal projection of $\R^n$ onto the 
span of the first $k$ standard basis vectors.  In this section,  
Theorem \ref{cont} is applied
to show that $P_k(Y)=(Y_,\ldots,Y_k)$ is approximately distributed as
a standard $k$-dimensional Gaussian random vector if $k=o(n)$. 
That $\E P_k(Y)P_k(Y)^T=I_k$ is immediate from the symmetry and
normalization, as above.  

This example is closely related to the
following result of Diaconis and Freedman in \cite{diafree}.  
\begin{thm}[Diaconis-Freedman]\label{diafreerankk}
Let $Z_1,\ldots,Z_n$ be independent standard Gaussian random variables and 
let $\P^k_\sigma$
be the law of $(\sigma Z_1,\ldots,\sigma Z_k)$.  
For a probability $\mu$ on $[0,\infty)$,
define $\P_{\mu,k}$ by
$$\P_{\mu,k}=\int\P^k_\sigma d\mu(\sigma).$$
Let $Y=(Y_1,\ldots,Y_n)\in\R^n$ be a spherically symmetric 
random vector, and let 
$\P_k$ be the law of $(Y_1,\ldots,Y_k)$.  Then
 there is
a probability measure $\mu$ on $[0,\infty)$ such that for $1\le k\le n-4,$
$$d_{TV}(\P_{k},\P_{\mu,k})\le\frac{2(k+3)}{n-k-3}.$$
Furthermore, the mixing measure $\mu$ can be taken to be the law of 
$\frac{1}{\sqrt{n}}|Y|.$
\end{thm}
In some cases, the explicit 
form given in Theorem \ref{diafreerankk}
for the mixing measure has allowed the theorem to 
be used to prove central limit theorems of 
interest in convex geometry; see \cite{bremv} and \cite{boaz}.
Theorem \ref{ksphere} below says that the variance bound \eqref{varbda} is 
sufficient to show that the mixing measure of Theorem \ref{diafreerankk}
can be taken to be a point mass.  In fact, it is not too 
difficult to obtain the total variation analog of Theorem \ref{ksphere} 
directly from the Diaconis-Freedman result and \eqref{varbda}; 
however, the Stein's method proof given 
below is considerably simpler than the direct proof given in 
\cite{diafree}.  The rates obtained are of the same order, though the rate 
obtained by 
Diaconis and Freedman is in the total variation distance, whereas  the rate
below is in the Wasserstein distance.  

\bigskip

To apply Theorem \ref{cont}, construct a family of exchangeable pairs 
as follows.  
For $\epsilon>0$ fixed, let 
\begin{equation*}\begin{split}
A_\epsilon&=\begin{bmatrix}\sqrt{1-\epsilon^2}&\epsilon\\-\epsilon&
\sqrt{1-\epsilon^2}\end{bmatrix}\oplus I_{n-2}\\
&=I_n+\begin{bmatrix}-\frac{\epsilon^2}{2}+\delta&\epsilon\\-\epsilon&
-\frac{\epsilon^2}{2}+\delta\end{bmatrix}\oplus0_{n-2},
\end{split}\end{equation*}
where $\delta$ is a deterministic constant and $\delta=O(\epsilon^4).$ 
Let $U$ be a Haar-distributed   $n\times n$ 
random orthogonal matrix, independent of $Y$, and let
$Y_\epsilon=\left(UA_\epsilon U^T\right)Y.$  Thus
$Y_\epsilon$ is a small random rotation of $Y$.  
In what follows, Theorem \ref{cont} is applied to the exchangeable pair
$(P_k(Y),P_k(Y_\epsilon)).$

Let $K$ be the $k\times 2$ matrix made of the first two 
columns of $U$ and
$C_2=\begin{bmatrix}0&1\\-1&0\end{bmatrix}.$  Define $Q:=KC_2K^T.$
Then by the construction of $Y_\epsilon,$
\begin{equation}\label{diff5}
P_k(Y_\epsilon)-P_k(Y)=\epsilon\left[-\left(\frac{\epsilon}{2}+\epsilon^{-1}
\delta\right)P_kKK^T+P_kQ\right]Y,
\end{equation}
and $\epsilon^{-1}\delta=O(\epsilon^3).$

To check the conditions of Theorem \ref{cont}, the following lemma is needed; 
see \cite{meckes-thesis}, Lemma 3.3 and Theorem 1.6 for a detailed proof.
\begin{lemma}\label{ints1}
If $U=\left[u_{ij}\right]_{i,j=1}^n$ is an 
orthogonal matrix distributed according to Haar measure, 
then $\E\left[\prod u_{ij}^{k_{ij}}\right]$ is non-zero
if and only if the number of entries from each row and from each
column is even.  Second and fourth-degree moments are as follows:
\begin{enumerate}
\item For all $i,j$, $$\E\left[u_{ij}^2\right]=\frac{1}{n}.$$ 
\item For all $i,j,r,s,\alpha,\beta,\lambda,\mu$,
\begin{equation*}\begin{split}
\E\big[u_{ij}u_{rs}&u_{\alpha\beta}u_{\lambda \mu}\big]\\&=
-\frac{1}{(n-1)n(n+2)}\Big[\delta_{ir}\delta_{\alpha\lambda}\delta_{j\beta}
\delta_{s\mu}+\delta_{ir}\delta_{\alpha\lambda}\delta_{j\mu}\delta_{s\beta}+
\delta_{i\alpha}\delta_{r\lambda}\delta_{js}\delta_{\beta\mu}\\&
\qquad\qquad\qquad\qquad\qquad\qquad+
\delta_{i\alpha}\delta_{r\lambda}\delta_{j\mu}\delta_{\beta s}+
\delta_{i\lambda}\delta_{r\alpha}\delta_{js}\delta_{\beta \mu}+
\delta_{i\lambda}\delta_{r\alpha}\delta_{j\beta}\delta_{s\mu}\Big]\\&\qquad
+\frac{n+1}{(n-1)n(n+2)}\Big[\delta_{ir}\delta_{\alpha\lambda}\delta_{js}
\delta_{\beta\mu}+\delta_{i\alpha}\delta_{r\lambda}\delta_{j\beta}\delta_{
s\mu}+\delta_{i\lambda}\delta_{r\alpha}\delta_{j\mu}\delta_{s\beta}
\Big].
\end{split}\end{equation*}
\label{hugemess}

\item For the matrix $Q=\big[q_{ij}\big]_{i,j=1}^n$ defined as above,
$q_{ij}=u_{i1}u_{j2}-u_{i2}u_{j1}.$
For all $i,j,\ell,p$,
$$\E\left[q_{ij}q_{\ell p}
\right]=\frac{2}{n(n-1)}\big[\delta_{i\ell}\delta_{jp}-\delta_{ip}\delta_{
j\ell}\big].$$
\end{enumerate}

\end{lemma}

By the lemma, $\E\big[KK^T\big]=\frac{2}{n}I$ and $\E\big[Q\big]=0,$
and so
\begin{equation*}\begin{split}
\lim_{\epsilon\to0}\frac{n}{\epsilon^2}\E\Big[(P_k(Y_\epsilon)-P_k(Y))\Big|
P_k(Y)\Big]&=-P_k(Y);
\end{split}\end{equation*}
condition \ref{first-diff} of Theorem \ref{cont} thus holds with $\lambda(
\epsilon)=\frac{\epsilon^2}{n}$.  

Fix $i,j\le k$.  By \eqref{diff5},
\begin{equation*}\begin{split}
\lim_{\epsilon\to0}\frac{n}{2\epsilon^2}\E\Big[(P_k(Y_\epsilon)-P_k(Y))_i&(
P_k(Y_\epsilon)-P_k(Y))_j\Big|Y\Big]\\&=\frac{n}{2}
\E\left[(P_kQY)_i(P_kQY)_j\big|Y\right]\\
&=\frac{n}{2}
\E\left.\left[\sum_{\ell, m}Y_\ell Y_mq_{i\ell}q_{jm}\right|Y\right]\\
&=\frac{1}{(n-1)}\E\left[\sum_{\ell, m}Y_\ell Y_m\left(\delta_{ij}
\delta_{\ell m}-\delta_{im}\delta_{\ell j}\right)\right]\\
&=\frac{1}{(n-1)}\left[\delta_{ij}|Y|^2 -Y_iY_j\right].
\end{split}\end{equation*}

Thus
$$F=\frac{1}{(n-1)}\left[\left(\E\left[|Y|^2-(n-1)\big|P_k(Y)\right]\right)
\cdot I_k-P_k(Y)P_k(Y)^T\right].$$
Now, 
$$\E\|P_k(Y)P_k(Y)^T\|_{H.S.}=\E\left|P_k(Y)\right|_2^2=k$$
by assumption, and $$\E\left|\E\left[|Y|^2-(n-1)\big|P_k(Y)\right]\right|\le \sqrt{a}+1,$$ so 
applying Theorem \ref{cont} gives:
\begin{thm}\label{ksphere}
With notation as above,  
$$d_{W}(P_k(Y),Z)\le\frac{k(\sqrt{a}+2)}{n-1}.$$
\end{thm}

\subsection{Rank $k$ projections of Haar measure on $\O_n$}

$\phantom{hi}$

\medskip

A theme in studying random matrices from the compact classical matrix groups
is that these matrices are in many ways (though not all ways) 
similar to Gaussian
random matrices.  For example, it was shown in 
\cite{limit} that if $M$ is a random matrix in the orthogonal group $\O_n$
distributed according to Haar measure, then
$$\sup_{\substack{A\,:\,\tr(AA^T)=n\\-\infty<x<\infty}}\big|\P(\tr(AM)\le x)-
\Phi(x)\big|\to0$$
as $n\to\infty.$  In \cite{M}, this result was refined to include a rate
of convergence (in total variation) of $W=\tr(AM)$ to a standard Gaussian 
random variable, depending only on the value of $\tr(AA^T)$.  That is, 
rank one projections of Haar measure on $\O_n$ are uniformly 
close to Gaussian, and rank one projections of Gaussian random matrices
are exactly Gaussian.

A natural question is whether rank $k$ projections of 
Haar measure on $\O_n$ are close, in some sense, to multivariate Gaussian
distributions, and if so, how large $k$ can be.  This is a more refined 
comparison of the type mentioned above, since the distributions of all 
projections of any rank of Gaussian matrices are Gaussian.  
In the remarkable recent 
work \cite{jiang06}, Tiefeng Jiang has shown that the entries of any 
$p_n\times q_n$ submatrix of an $n\times n$ random orthogonal matrix  are 
close to i.i.d.\ Gaussians in total variation distance whenever $p_n = 
o(\sqrt{n})$ and $q_n  = o(\sqrt{n})$, and that these orders of $p_n$ and $q_n$
are best possible.  This improved an earlier result of 
Diaconis, Eaton, and Lauritson \cite{del92}, which proved the result in the 
case of $p_n = o(n^{1/3})$ and $q_n  = o(n^{1/3})$.  
As this article was in preparation, Beno\^it Collins and Michael Stolz
\cite{cs06} proved that for $r$ fixed, $A_1,\ldots,A_r$ deterministic 
parameter matrices,
and $M$ a uniformly distributed element of a classical compact
symmetric space (represented as a space of matrices), the random vector 
$(\tr(A_1M),\ldots,\tr(A_rM))$ converges weakly to a Gaussian random vector,
as the dimension of the space tends to infinity.  Their work in particular
covers the cases of $M$ a Haar-distributed random orthogonal or unitary matrix,
but goes farther to consider more general homogeneous spaces.

In this section, it is shown
that rank $k$ projections of Haar measure on $\O_n$ are close in 
Wasserstein distance to Gaussian
for $k=o(n)$.    This in particular recovers Jiang's 
result (in Wasserstein distance), but is more general in that it is uniform
over all rank $k$ projections, and not just those having the special form of 
truncation to a sub-matrix.  
The theorem also strengthens 
the result of Collins and Stolz, in the case that $M$ is 
a random element of $\O_n$.

\begin{thm}\label{mix}
Let $B_1,\ldots,B_k$ be linearly independent $n\times n$ matrices (i.e.\ the only linear combination of them 
which is equal to the zero matrix has all coefficients equal to zero) 
over $\R$ such that 
$\tr(B_iB_i^T)=n$ for each $i$.  Let $b_{ij}=\tr(B_iB_j^T).$  Let $M$ be
a random orthogonal matrix and let $$X=(\tr(B_1M),\tr(B_2M),\ldots, \tr(B_kM))
\in\R^k.$$  Let $Y=(Y_1,\ldots,Y_k)$ be a random vector whose components have
the 
standard Gaussian distribution, with covariance matrix 
$C:=\frac{1}{n}\left(b_{ij}
\right)_{i,j=1}^k$.  Then for $n\ge 2,$ 
$$d_{W}(X,Y)\le\frac{k\sqrt{2\|C\|_{op}}}{n-1}.$$
\end{thm}

\medskip

\noindent{\em Remark.}
Lemma \ref{ints1} and an easy computation show that 
for all $i,j$,
$$
\E\big[\tr(B_iM)\tr(B_jM)\big]=\frac{1}{n}\inprod{B_i}{B_j},$$
thus the matrix $C$ above is also the covariance matrix of $X$.

\bigskip

It is shown below that Theorem \ref{mix} follows fairly easily from the 
following special case.
\begin{thm}\label{ind}
Let $A_1,\ldots,A_k$ be $n\times n$ matrices over $\R$ satisfying
$\tr(A_iA_j^T)=n\delta_{ij};$ for $i\neq j$, $A_i$ and $A_j$ are 
orthogonal with respect to the Hilbert-Schmidt inner product.
Let $M$ be a random orthogonal matrix, and consider the vector
$X=(\tr(A_1M),\tr(A_2M),\ldots, \tr(A_kM))\in\R^k.$  
Let $Z=(Z_1,\ldots,Z_k)$ be a random vector whose components are independent
standard normal random variables.  Then for $n\ge 2$,

$$\big|\E f(X)-\E f(Z)\big|\le\frac{\sqrt{2}M_1(f)k}{n-1}$$
where $M_1(f)$ is the Lipschitz constant of $f$.  
\end{thm}
\subsection*{Example}
Let $M$ be a random $n\times n$ orthogonal matrix, and let 
$0<a_1<a_2<\ldots<a_k=n$.  For each $1\le i\le n,$ let 
$$B_i=\sqrt{\frac{n}{a_i}}I_{a_i}\oplus{\bf 0}_{n-a_i};$$
$B_i$ has $\sqrt{\frac{n}{a_i}}$ in the first $a_i$ diagonal entries and
zeros everywhere else.  If $i\le j$, then $\inprod{B_i}{B_j}_{HS}
=n\sqrt{\frac{a_i}{a_j}};$ in particular, $\inprod{B_i}{B_i}_{HS}=n$.
The $B_i$ are linearly independent w.r.t. the Hilbert-Schmidt inner
product since the $a_i$ are all distinct, so to apply Theorem \ref{mix},
we have only to bound the eigenvalues of the matrix $\left(\sqrt{\frac{a_{
\min(i,j)}}{a_{\max(i,j)}}}\right)_{i,j=1}^k.$  But this is easy, since
$|\lambda|\le\sqrt{\sum_{i,j=1}^k\frac{a_{\min(i,j)}}{a_{\max(i,j)}}}\le k$
for all eigenvalues $\lambda$ (see, e.g., \cite{bhat}).  It now follows
from Theorem \ref{mix} that if $Y$ is a vector of standard normals with
covariance matrix $\left(\sqrt{\frac{a_{\min(i,j)}}{a_{\max(i,j)}}}
\right)_{i,j=1}^k$ and $X=(\tr(B_1M),\ldots,\tr(B_kM))$, then 
$$\sup_{|f|_L\le1}\big|\E f(X)-\E f(Y)\big|\le\frac{\sqrt{2}k^{3/2}
}{n-1}.$$

\section*{Proofs}
\begin{proof}[Proof of Theorem \ref{mix} from Theorem \ref{ind}]
Perform the Gram-Schmidt algorithm on the matrices $\{B_1,\ldots,B_k\}$
with respect to the Hilbert-Schmidt inner product $\inprod{C}{D}=\tr(CD^T)$
to get matrices $\{A_1,\ldots,A_k\}$ which are mutually orthogonal and 
have H-S norm $\sqrt{n}$.  Denote the matrix which takes
the $B$'s to the $A$'s by 
$D^{-1}$ for $D=\big[d_{ij}\big]_{i,j=1}^n$; 
the matrix is invertible since the 
$B$'s are linearly independent.  Now by assumption,
\begin{equation*}\begin{split}
b_{ij}&=\inprod{B_i}{B_j}\\
&=\inprod{\sum_ld_{il}A_l}{\sum_pd_{jp}A_p}\\
&=n\sum_{l}d_{il}d_{jl}.
\end{split}\end{equation*}
Thus $DD^T=C=\frac{1}{n}\left(b_{ij}\right)_{i,j=1}^k.$

Now, let $f:\R^k\to\R$ with $M_1(f)\le1$.  Define $h:\R^k\to\R$ by 
$h(x)=f(Dx).$  Then $M_1(h)\le\|D\|_{op}\le\sqrt{\|DD^T\|_{op}}$.
By Theorem \ref{ind},
$$\big|\E h(\tr(A_1M),\ldots,\tr(A_kM))-\E h(Z)\big|\le\frac{k
\sqrt{2\|C\|_{op}}}{n-1}$$
for $Z$ a standard Gaussian random vector in $\R^k$.  But 
$D\big(\tr(A_1M),\ldots,\tr(A_kM)\big)=\big(\tr(B_1M),\ldots,\tr(B_kM)\big)$
and $DZ$ has standard normal components with covariance matrix 
$C=\frac{1}{n}\left(b_{ij}\right)_{i,j=1}^k.$
\end{proof}

\bigskip

\begin{proof}[Proof of Theorem \ref{ind}]

Make an exchangeable
pair $(M,M_\epsilon)$ as before; let $A_\epsilon$ be the rotation
$$A_\epsilon=
\begin{bmatrix}\sqrt{1-\epsilon^2}&\epsilon\\-\epsilon&\sqrt{1-\epsilon^2}
\end{bmatrix}\oplus I_{n-2}=I_n+
\begin{bmatrix}
\sqrt{1-\epsilon^2}-1&\epsilon\\-\epsilon&\sqrt{1-\epsilon^2}-1\end{bmatrix}
\oplus{\bf 0}_{n-2},$$
let $U$ be a Haar-distributed 
random orthogonal matrix, 
independent of $M$, and let 
$$M_\epsilon=UAU^TM.$$
Let $X_\epsilon=(\tr(A_1M_\epsilon),\ldots, \tr(A_kM_\epsilon)).$

As in section \ref{sphere}, define $K$ to be the first two columns of $U$ 
and $C_2=\begin{bmatrix}0&1\\-1&0
\end{bmatrix}$, and let $Q=KC_2K^T.$  Then
\begin{equation}\label{diff4}
M_\epsilon-M=\epsilon\left[\left(\frac{-\epsilon}{2}+O(\epsilon^3)\right)
KK^T+Q\right]M.\end{equation}
It follows from Lemma \ref{ints1} that
$\E\big[KK^T\big]=\frac{2}{n}I$ and $\E\big[Q\big]=0,$
thus
\begin{equation*}\begin{split}\lim_{\epsilon\to0}\frac{n}{\epsilon^2}
\E\big[(X_\epsilon&-X)_i\big|M\big]\\&=\lim_{\epsilon
\to0}\frac{n}{\epsilon^2}
\E\left[\tr[A_i(M_\epsilon-M)]\big|M\right]\\&=\lim_{\epsilon\to0}
\frac{n}{\epsilon^2}\left[\left(-\frac{\epsilon^2}{2}+O(\epsilon^4)\right)
\E\left[\tr(A_iKK^TM)\big|M\right]+\epsilon\E\left[\tr(A_iQM)\big|M
\right]\right]\\&=\lim_{\epsilon\to0}\frac{n}{\epsilon^2}
\left(-\frac{\epsilon^2}{2}+O(\epsilon^4)\right)\frac{2}{n}X_i\\
&=-X_i.\end{split}\end{equation*}
Condition \ref{first-diff} of Theorem \ref{cont} is thus satisfied with $\lambda(\epsilon)=\frac{\epsilon^2}{n}$.
The random matrix $F$ is computed as follows.  For notational
convenience, write $A_i=A=(a_{pq})$ and $A_j=B=(b_{\alpha\beta}).$  
By (\ref{diff4}), 
\begin{equation}\begin{split}\label{ij}
\lim_{\epsilon\to0}\frac{n}{2\epsilon^2}\E\Big[(X_\epsilon-X)_i&
(X_\epsilon-X)_j\Big|M\Big]\\&=\frac{n}{2}\E\left[\tr(AQM)\tr(BQM)\big|M
\right]\\&=\frac{n}{2}\E\left[\left.\sum_{p,q,r,\alpha,\beta,\gamma}a_{pq}
b_{\alpha\beta}m_{rp}m_{\gamma\alpha}q_{qr}q_{\beta\gamma}\right|M\right]\\&=
\frac{n}{2}\E\left[\sum_{p,q,r,\alpha,\beta\gamma}a_{pq}b_{\alpha\beta}
m_{rp}m_{\gamma\alpha}\left(\frac{2}{n(n-1)}\right)(\delta_{q\beta}\delta_{r
\gamma}-\delta_{q\gamma}\delta_{r\beta})\right]\\&=\frac{1}{(n-1)}\E\left[\inprod{MA}{MB}_{H.S.}-
\tr(AMBM)\right]\\&=\frac{1}{(n-1)}
\E\left[\inprod{A}{B}_{H.S.}-\tr(MAMB)\right]\\&=\frac{1}{(n-1)}
\left[n\delta_{ij}-\tr(MAMB)\right].
\end{split}\end{equation} 
Thus
$$F=\frac{1}{(n-1)}\E\left[\Big[\delta_{ij}-\tr(A_iMA_jM)
\Big]_{i,j=1}^k\Big|X\right].$$

\bigskip

{\bf Claim:} 
If $n\ge2 $, then $\E\left[\tr(A_iMA_jM)-\delta_{ij}\right]^2
\le 2$ for all $i$ and $j$.   

\bigskip

The claim gives that, for $n\ge 2$,
$$\E\|F\|_{H.S.}\le\sqrt{\E\|F\|_{H.S.}^2}\le \frac{\sqrt{2}k}{n-1},$$
thus completing the proof.

To prove the claim, 
first observe that Lemma \ref{ints1} implies 
$$\E\big[\tr(A_iMA_jM)\big]=\frac{1}{n}\inprod{A_i}{A_j}=\delta_{ij}.$$
Again writing $A_i=A$ and $A_j=B$, applying Lemma \ref{ints1} gives, 
\ref{hugemess},
\begin{eqnarray*}
\E\left[\tr(AMBM)\right]^2&=&\E\left[\sum_{\substack{p,q,r,s\\\alpha,\beta,
\mu,\lambda}}a_{sp}a_{\mu\alpha}b_{qr}b_{\beta\lambda}m_{pq}m_{rs}m_{\alpha
\beta}m_{\lambda\mu}\right]\\&=&-\frac{2}{(n-1)n(n+2)}\left[\tr(A^TAB^TB)+
\tr(AB^TAB^T)+\tr(AA^TBB^T)\right]\\&&\qquad
+\frac{n+1}{(n-1)n(n+2)}\left[2\inprod{A}
{B}_{H.S.}+\|A\|_{H.S.}^2\|B\|_{H.S.}^2\right]
\end{eqnarray*}
Now, as the Hilbert-Schmidt norm is submultiplicative (see \cite{bhat}, page 
94),
$$\tr(A^TAB^TB)\le\|A^TA\|_{H.S.}\|B^TB\|_{H.S.}\le\|A\|_{H.S.}^2
\|B\|_{H.S.}^2=n^2,$$
and the other two summands of the first line are
bounded by $n^2$ in the same way.  Also,
$$2\inprod{A}{B}_{H.S.}+\|A\|_{H.S.}^2\|B\|_{H.S.}^2=n^2(1+2\delta_{ij}),$$
Thus 
$$\E\left[\tr(A_iMA_jM)-\delta_{ij}\right]^2\le
\frac{-6n^2+(n+1)n^2(1+2\delta_{ij})-(n-1)n(n+2)\delta_{ij}}{(n-1)n(n+2)}
\le 2.$$
\end{proof}

\subsection{Complex-linear functions of random unitary matrices}\label{unitary}
$\phantom{hi}$

\medskip

In this section, we consider Haar-distributed random matrices in $\U_n$.
As discussed in the previous section, a general theme in studying random
matrices from the classical compact matrix groups has been to compare
to the corresponding Gaussian distribution.   In particular, 
it was shown in \cite{limit} that if 
$M=\Gamma+i\Lambda$ is a random $n\times n$ unitary matrix and 
$A$ and $B$ are fixed real diagonal matrices with $\tr (AA^T)=\tr(BB^T)=n,$
then $\tr(A\Gamma)+i\tr(B\Lambda)$ converges in distribution to standard
complex normal.  This implies in particular that $\Re(\tr(AM))$ converges in
distribution to $\n\left(0,\frac{1}{2}\right).$  A total variation 
rate of convergence for
this last statement was obtained in \cite{M}, giving as an easy consequence
the weak-star convergence of the random variable $W=\tr(AM)$ to standard
complex normal, for $A$ an 
$n\times n$ matrix over $\C$ with $\tr(AA^*)=n$.  The approaches used in
\cite{limit} and \cite{M} are somewhat awkward, partly due to the fact that
the limiting behavior of $W$ is a multivariate question.  In
this section, Corollary \ref{complex} is applied to prove the analogous 
result to Theorem \ref{ind} for complex-rank $k$ projections of 
Haar measure on the space of random unitary matrices.  
As in the previous section, this result recovers and strengthens the result
of Collins and Stolz \cite{cs06}, in the case that $M$ is a Haar-distributed
unitary matrix.

\medskip

\begin{thm}\label{uthm}
Let $M\in\U_n$ be distributed according to Haar measure, and let 
$\{A_i\}_{i=1}^k$ be fixed $n\times n$ matrices over $\C$ such that 
$\tr(A_iA_j^*)=n\delta_{ij}$.  Let $W(M)=(\tr(A_1M),\ldots,\tr(A_kM))$ and let
$Z$ be a standard complex Gaussian random vector in $\C^k$. 
Then there is a universal constant $c$ such that 
\begin{equation*}
d_{W}(W,Z)\le\frac{ck}{n}.
\end{equation*}

\end{thm}

\medskip

{\em Remark:}  The constant $c$ given by the proof  
is asymptotically
equal to $\sqrt{2}$; for $n\ge 4$, $c$ can be taken to be 3.

\medskip

For the proof, the following lemma is needed.  See \cite{meckes-thesis}, Lemma
3.5 for a detailed proof.
\begin{lemma}\label{unit-ints}
Let $H=\big[h_{ij}\big]_{i,j}\in\U_n$ be distributed according to Haar measure.  Then
the expected value of a product of entries of $H$ and their conjugates
is non-zero only when there are the same number of entries as 
conjugates of entries from each row and from each column.  
Second- and fourth-degree moments are as follows.
\begin{enumerate}
\item For all $i,j$, $$\E\big[|h_{ij}|^2\big]=\frac{1}{n},$$ \label{norm2}
\item For all $i,j,r,s,\alpha,\beta,\lambda,\mu$,\label{degree4}
\begin{equation*}\begin{split}
\E\big[h_{ij}h_{rs}\overline{h}_{\alpha \beta}\overline{h}_{\lambda \mu}\big]&=
\frac{1}{(n-1)(n+1)}\Big[\delta_{i\alpha}\delta_{r\lambda}\delta_{j\beta}\delta_{
s\mu}+\delta_{i\lambda}\delta_{r\alpha}\delta_{j\mu}\delta_{s\beta}\Big]\\&\qquad
\qquad-\frac{1}{(n-1)n(n+1)}\Big[\delta_{i\alpha}\delta_{r\lambda}\delta_{j\mu}
\delta_{s\beta}+\delta_{i\lambda}\delta_{r\alpha}\delta_{j\beta}\delta_{s\mu}
\Big],
\end{split}\end{equation*}
\item \label{twist-mean}
\begin{equation*}\begin{split}
\E\big[(h_{i1}\overline{h}_{j2}-h_{i2}\overline{h}_{j1})
(h_{r1}\overline{h}_{s2}-h_{r2}\overline{h}_{s1})\big]&=
-\frac{2}{(n-1)(n+1)}\,\delta_{is}\delta_{jr}+\frac{2}{(n-1)n(n+1)}\,\delta_{ij}
\delta_{rs}.\end{split}\end{equation*} 
\end{enumerate}
\end{lemma}

\begin{proof}[Proof of Theorem \ref{uthm}]

The theorem is proved as an application of Corollary \ref{complex}, similarly
to the proof of Theorem \ref{ind} via Theorem \ref{cont}.  Construct a 
family of pairs
$(W,W_\epsilon)$ analogously to what was done in the orthogonal case:
let $U\in\U_n$ be a random unitary matrix, independent of 
$M$, and let $M_\epsilon=UA_\epsilon U^*M$, where as before
$$A_\epsilon=\begin{bmatrix}\sqrt{1-\epsilon^2}&\epsilon\\-\epsilon&\sqrt{1-
\epsilon^2}\end{bmatrix}\oplus I_{n-2},$$
thus $M_\epsilon$ is a small random rotation of $M$.  Let $W_\epsilon=W(M_
\epsilon)$; $(W,W_\epsilon)$ is exchangeable by construction.

As in the previous sections, let
$I_2$ be the $2\times2$ identity matrix, $K$ the $n\times
2$ matrix made from the first two columns of $U=\big[u_{ij}\big]_{i,j}$, and let 
$$C_2=\begin{bmatrix}0&1\\-1&0\end{bmatrix}.$$  
Define the matrix $Q=KC_2K^*$.
Then 
\begin{equation*}
M_\epsilon=M+K\big[(\sqrt{1-\epsilon^2}-1)I_2+\epsilon C_2\big]K^*M,
\end{equation*}
and
\begin{equation}\label{diff3}
\tr(A_iM_\epsilon)-\tr(A_iM)=\left(-\frac{\epsilon^2}{2}+O(\epsilon^4)\right)
\tr(A_iKK^*M)+\epsilon\tr(A_iQM).
\end{equation}
It follows from Lemma \ref{unit-ints} that 
$\E \big[KK^*\big]=\frac{2}{n}I$
and
$\E \big[Q\big] =0,$  
thus
\begin{equation}\label{grad}
\begin{split}
\lim_{\epsilon\to0}\frac{n}{\epsilon^2}\E\left[\tr(A_iM_\epsilon)-\tr(A_iM)
\big|M\right]&=-\tr(A_iM),
\end{split}\end{equation}
and the first condition of Corollary \ref{complex} holds with $\lambda(\epsilon)
=\frac{\epsilon^2}{n}$.

Let $A_i=:A=\big[a_{pq}\big]$ and $A_j=:B=\big[b_{\alpha\beta}\big]$; 
by (\ref{diff3}) and Lemma \ref{unit-ints},
\begin{equation}\begin{split}
\lim_{\epsilon\to0}\frac{n}{2\epsilon^2}\E&\left[(W_\epsilon 
-W)_i(W_\epsilon-W)_j\big|W\right]\\&=\frac{n}{2}\E\left[\left(\tr(AQM)\right)\left(\tr(BQM)
\right)\big|W\right]\\
&=\frac{n}{2}\E\left.\left[\sum_{p,q,r,\alpha,\beta,\mu}a_{pq}m_{rp}b_{\alpha\beta}m_{\gamma
\alpha}(u_{q1}
\overline{u}_{r2}-u_{q2}\overline{u}_{r1})(u_{\beta1}
\overline{u}_{\gamma 2}-u_{\beta2}\overline{u}_{\gamma1})\right|W\right]\\&=
\frac{1}{(n-1)(n+1)}\left[\sum_{p,q,\alpha,\beta}a_{pq}m_{qp}b_{\alpha
\beta}m_{\beta\alpha}-n\sum_{p,q,\alpha\beta}a_{pq}b_{\alpha\beta}m_{\beta p}m_{q\alpha}\right]\\&=\frac{1}{(n-1)(n+1)}\left[
\tr(AM)\tr(BM)-n\tr(AMBM)\right].\label{secdiff1}
\end{split}\end{equation}
Similarly, one can use part \ref{twist-mean} of Lemma \ref{unit-ints} with 
the roles of $r$ and $s$ reversed to get
\begin{equation}\begin{split}\label{mnormdiff1}
\lim_{\epsilon\to0}\frac{n}{2\epsilon^2}\E\big[(W_\epsilon& -W)_i\overline{
(W_\epsilon-W)}_j\big|W\big]=\frac{n}{2}
\E\left[\tr(AQM)\overline{\tr(BQM)}\big|W\right]\\
&=\frac{n}{2}\E\left.\left[\sum_{p,q,r,\alpha,\beta,\gamma}a_{pq}m_{rp}
\overline{b}_{\alpha\beta}\overline{m}_{\gamma\alpha}(u_{q1}
\overline{u}_{r2}-u_{q2}\overline{u}_{r1})(\overline{u}_{\beta 1}
u_{\gamma2}-\overline{u}_{\beta2}u_{\gamma1})\right|W\right]\\&=
\frac{1}{(n-1)(n+1)}\left[n\sum_{p,\alpha}\left(\sum_qa_{pq}\overline{b}_{\alpha 
q}\right)\left(\sum_\gamma m_{\gamma p}\overline{m}_{\gamma\alpha}\right)-
\sum_{p,q,\alpha,\beta}a_{pq}\overline{b}_{\alpha\beta}m_{qp}
\overline{m}_{\beta\alpha}\right]\\&=
\frac{1}{(n-1)(n+1)}\left[n^2\delta_{ij}-\tr(AM)\overline{\tr(BM)}
\right]\\&=\delta_{ij}
+\frac{1}{(n-1)(n+1)}\Big[\delta_{ij}-\tr(AM)\overline{\tr(BM)}\Big],
\end{split}\end{equation}
where the fact that $M$ is unitary and the assumption $\tr(A_iA_j^*)=n
\delta_{ij}$ have been used to get the second to last line.

One can thus take
\begin{equation*}\begin{split}
\gamma_{ij}=\frac{\delta_{ij}-\tr(A_iM)\overline{\tr(A_jM)}}{(n-1)(n+1)}
\qquad\quad\lambda_{ij}=\frac{\tr(A_iM)\tr(A_jM)-n\tr(A_iMA_jM)}{(n-1)(n+1)}.
\end{split}\end{equation*}

By the Cauchy-Schwarz inequality, 
$$\E\|\Gamma\|_{H.S.}\le\sqrt{\E\sum_{i,j}|\gamma_{ij}|^2}$$
and
\begin{equation*}\begin{split}
\E|\gamma_{ij}|^2&=\frac{1}{(n-1)^2(n+1)^2}\E\left[\delta_{ij}-2\Re(\tr(A_iM)
\overline{\tr(A_jM)})+|\tr(A_iM)\tr(A_jM)|^2\right].
\end{split}\end{equation*}
Now,
$$\E|\tr(A_iM)\overline{\tr(A_jM)}|\le\sqrt{\E|\tr(A_iM)|^2\E|\tr(A_jM)|^2}=1$$
by the normalization of the matrices $A_i$.  Again writing $A=A_i$ and $B=A_j$,
\begin{equation*}\begin{split}
\E|&\tr(AM)\tr(BM)|^2\\&=\sum_{\substack{p,q,r,s\\\alpha,\beta,\lambda,\mu}}
a_{pq}\overline{a}_{rs}b_{\alpha\beta}\overline{b}_{\lambda\mu}\E\big[m_{qp}m_{
\beta\alpha}\overline{m}_{sr}\overline{m}_{\mu\lambda}\big]\\&=\frac{1}{(n-1)
(n+1)}\left[\tr(AA^*)\tr(BB^*)+(\tr(AB^*))^2-\frac{1}{n}\tr(AA^*BB^*)-\frac{1}
{n}\tr(A^*AB^*B)\right]\\&=\frac{1}{(n-1)(n+1)}\left[n^2(1+\delta_{ij})-\frac{1}
{n}\tr(AA^*BB^*)-\frac{1}{n}\tr(A^*AB^*B)\right],
\end{split}\end{equation*}
where Lemma \ref{unit-ints} has been used to get the third line and the
normalization and orthogonality conditions on the $A_i$ have been used to 
get the last line.  Now,
$$|\tr(AA^*BB^*)|\le\|AA^*\|_{H.S.}\|BB^*\|_{H.S.}\le\|A\|_{H.S.}\|A^*\|_{H.S.}
\|B\|_{H.S.}\|B^*\|_{H.S.}=n^2;$$
the first inequality is just the Cauchy-Schwarz inequality for the 
Hilbert-Schmidt inner product and the second is due to the submultiplicativity
of the Hilbert-Schmidt norm (see \cite{bhat}, page 94).  It now follows that
$$\E|\gamma_{ij}|^2\le\frac{1}{(n-1)^2(n+1)^2}\left[\delta_{ij}+2+\frac{n^2(1+
\delta_{ij})+2n}{(n-1)(n+1)}\right]\le\frac{1}{(n-1)^2(n+1)^2}\left[5+\frac{2}{
n-1}\right],$$
and thus
\begin{equation}\label{gamma-bd}
\E\|\Gamma\|_{H.S.}\le\frac{k}{(n-1)(n+1)}\sqrt{5+\frac{2}{n-1}}.
\end{equation}

Taking a similar approach to bounding $\E\|\Lambda\|_{H.S.}$,
\begin{equation}\begin{split}\label{lambda-1}
\E|\lambda_{ij}|^2=\frac{1}{(n-1)^2(n+1)^2}
\E\Big[|\tr(A_iM)\tr(A_jM)|^2&-2n\Re(\tr(A_iM)\tr(A_jM)\overline{\tr(A_iMA_jM)}
)\\&+n^2|\tr(A_iMA_jM)|^2\Big].
\end{split}\end{equation}
It has already been shown that 
$$\E|\tr(A_iM)\tr(A_jM)|^2\le\frac{n^2(1+\delta_{ij})+2n}{(n-1)(n+1)}\le
2+\frac{2}{n-1}.$$
One can use Lemma \ref{unit-ints} to compute the other two terms similarly:
\begin{equation*}\begin{split}
\E&\Big[\tr(AM)\tr(BM)\overline{\tr(AMBM)}\Big]\\
&=\sum_{\substack{p,q,r,s\\\alpha,
\beta, \lambda, \mu}}a_{pq}\overline{a}_{\alpha\beta}b_{rs}\overline{b}_{\lambda
\mu}\E\big[m_{qp}m_{sr}\overline{m}_{\beta\lambda}\overline{m}_{\mu\alpha}\big]\\
&=\frac{1}{(n-1)(n+1)}\left[\tr(AA^*BB^*)+\tr(A^*AB^*B)-\frac{1}{n}\tr(AA^*)
\tr(BB^*)-\frac{1}{n}(\tr(AB^*))^2\right]\\&=\frac{1}{(n-1)(n+1)}\left[
\tr(AA^*BB^*)+\tr(A^*AB^*B)-n(1+\delta_{ij})\right],
\end{split}\end{equation*}
thus
$$\left|\E\left[2n\Re\Big(\tr(A_iM)\tr(A_jM)\tr(A_iMA_jM)\Big)\right]\right|
\le\frac{4n^3+2n(1+\delta_{ij})}{(n-1)(n+1)}\le4n+\frac{8}{n-1};$$
and
\begin{equation*}\begin{split}
\E&|\tr(AMBM)|^2\\&=\sum_{\substack{p,q,r,s\\\alpha,\beta,\lambda,\mu}}
a_{pq}\overline{a}_{\alpha\beta}b_{rs}\overline{b}_{\lambda\mu}\E\big[m_{qr}m_{sp}
\overline{m}_{\beta\lambda}\overline{m}_{\mu\alpha}\big]\\&=\frac{1}{(n-1)(n+1)}
\left[\tr(AA^*)\tr(BB^*)+(\tr(AB^*))^2-\frac{1}{n}\tr(AA^*BB^*)-\frac{1}{n}
\tr(A^*AB^*B)\right]\\&=\frac{1}{(n-1)(n+1)}
\left[n^2(1+\delta_{ij})-\frac{1}{n}\tr(AA^*BB^*)-\frac{1}{n}
\tr(A^*AB^*B)\right],
\end{split}\end{equation*}
thus
$$n^2\E|\tr(A_iMA_jM)|^2\le\frac{n^4(1+\delta_{ij})+2n^3}{(n-1)(n+1)}\le
2n^2+\frac{2n^2}{n-1}.$$
Using these three bounds in \eqref{lambda-1} yields
$$\E\|\Lambda\|_{H.S.}\le\sqrt{\sum_{i,j}\E|\lambda_{ij}|^2}\le
\frac{k}{(n-1)}\sqrt{2+\frac{2(n^2+5)}{(n-1)(n+1)^2}}.$$
\end{proof}

\bibliographystyle{plain}

\end{document}